\documentclass[10pt]{amsart}
\usepackage{epsfig,amsfonts,amsthm,amssymb,latexsym,amsmath,enumerate,hyperref,color}

\usepackage{epsfig}
\usepackage{epstopdf}

\newtheorem{theorem}{Theorem}
\newtheorem{corollary}[theorem]{Corollary}

\newtheorem{remark}{Remark}

\title[Effective Methods on Systems of Nonlinear Difference Equations]
{Effective Methods on Determining the Periodicity and Form of Solutions of Some Systems of Nonlinear Difference Equations}

\author[J. F. T. Rabago]{Julius Fergy T. Rabago}
\address{Julius Fergy Tiongson Rabago, Department of Mathematics and Computer Science,
College of Science, University of the Philippines Baguio, Baguio City 2600, Benguet, PHILIPPINES} \email{jfrabago@gmail.com}

\subjclass[2010]{Primary: 39 A 10; Secondary: 37 B 20}

\keywords{Difference equations, periodicity of solutions, system of equations.\\ 
\indent {\it Preprint submitted for publication on} November 18, 2015.}

\date{November 18, 2015}
\begin{document}
\maketitle

\begin{abstract}
Recently, various systems of nonlinear difference equations, of different forms, were studied.
In this existing work, two earlier published papers, due respectively to Bayram and Da\d{s} [Appl. Math. Sci. (Ruse), {\bf 4}(7) (2010) pp. 817--821] and Elsayed [Fasciculi Mathematici, {\bf 40} (2008), pp. 5--13], are revisited. 
The results exhibited in these previous investigations are re-examined through a new approach, more theoretical and explanative compared to the ones offered in these aforementioned works.
Furthermore, the qualitative behavior of solutions of a system of nonlinear difference equations of higher-order is investigated through analytical methods.
The system, which is considered here, generalizes those that are first presented in [Fasciculi Mathematici, {\bf 40} (2008), pp. 5--13] but are treated differently from this pre-existing work. 
The results delivered here are important, not only for the reason that they provide theoretical explanations on several earlier results, but also because the system being studied can be use to model real-life phenomena exhibiting periodic behaviors. 

\end{abstract}


\section{Introduction}
Sequences often arise from recursion relations. 
Difference equations, which have been of great interest in previous decades, are known to fall in this category.
In fact, this type of equation has attracted the attentions of many mathematicians in recent years. 
This curiosity, as we have witnessed, continuous to grow as more fascinating results are obtained and presented in recent studies.
Interestingly, in some cases, difference equations may be solve in closed-forms (see, e.g., \cite{rabago1, touafek, tollu, tollu1, tollu2}).
Moreover, in some instances, the solutions of these equations may be shown to be periodic.
In regards to this intriguing behavior, we say that a solution of a given difference equation, say $x_{n+1}=f(x_,x_{n-1},\ldots,x_{n-k})$ (which is of order $k+1$), 
is \emph{periodic} or \emph{eventually periodic} with period $p>0$ if we can find an integer $N\geq -k$
such that the relation $x_{n+p}=x_n$ holds for all $n\geq N$ (cf. \cite{ladas}).

Difference equations may be classified as \emph{linear} and \emph{nonlinear} types.
A classical example of a linear type is the well-known Fibonacci sequence $\{F_n\}:=\{F_n\}_{n=0}^{\infty}=\{0,1,1,2,3,5,7,\ldots\}$ (cf. \cite{koshy})
which has been generalized by many authors in various forms (see, for instance, the so-called Horadam sequence \cite{horadam,lbf,lucas}). 
Quite recently, these sequences have been applied in the study of functional and differential equations (see \cite{rabago1, rabago2, rabago3}, and the references therein).
On the other hand, nonlinear types may be regarded as ratios of linear difference equations.
The simplest example of such type is the difference equation 
\begin{equation}
\label{ex}
x_{n+1} =\frac{1}{x_n}, \qquad n=0,1,\ldots,
\end{equation}
with a real initial value $x_0 \neq 0$.
We mention that the $n$-th term of the Fibonacci sequence $\{F_n\}$ is known to be solvable explictly in $n$, 
i.e., we may express the $n$-th term Fibonacci number $F_n$ in Binet's form \cite{weisstein}:
\[
F_n = \frac{\phi^n-(-\phi)^{-n}}{\sqrt{5}}, \qquad n=0,1,\ldots,
 \] 
where $\phi=(1+\sqrt{5})/2=1.61803398874\ldots$ denotes the widely known \emph{golden ratio}.
Clearly, the Fibonacci sequence is unbounded and grows exponentially to infinity as $n$ increases without bound.
In contrast to this behavior, every solutions of the difference equation \eqref{ex}, regardless of its initial value, is periodic with period $2$.
This is seen easily as follows: since $x_{n+1}=1/x_n$, then iterating the left hand side one-step lower yields the relation $x_{n+1}=x_{n-1}$, from which it is evident that the sequence is periodic with periodicity $2$.
Knowing if a sequence is periodic or not in advance has some sort of advantage in solving systems of difference equations.
It will, possibly, lessen ones effort in establishing the form of solutions of these types of equations.
For example, if we know for instance that a solution of a difference equation is eventually periodic, starting from some index $n_0$, with period $p$ and we find it difficult to establish its solution form, 
then, perhaps, we can at least manually find their forms by iteratively computing term by term their respective forms up to the $n_0+p$-th term of the solution sequence.

Arguably, difference equations have many applications, not only in various fields of mathematics, but also in related sciences (see, e.g., a very recent book by Mickens \cite{mickens}).
In fact, as it is popularly known, the origin of Fibonacci sequence outside India \cite{knuth} came from the book \emph{Liber Abaci} (1202) by Leonardo of Pisa 
from which the sequence is describe as the ideal (biologically unrealistic) growth pattern of rabbits (cf. \cite{koshy}).  
Moreover, the sequence which is known to be generated by the difference equation $F_{n+1}=F_n+F_{n-1}$ (with initial values $F_0=F_1=1$) has applications in computer algorithm (cf. \cite{knuth}). 
They also appear in biological settings, such as branching in trees, {\it phyllotaxis} or the arrangement of leaves of stem \cite{douady},
the fruit sprouts of a pineapple \cite{jones}, the flowering of an artichoke, an uncurling fern and the arrangement of a pine cone's bracts \cite{brousseau}, etc.
Nonlinear types of difference equations also posses numerous applications in pure mathematics and related sciences.
For instance, the Newton's method (or Newton-Raphson's method), which is known as a root-finding algorithm, uses the nonlinear difference equation
\[
x_{n+1}=x_n-\frac{f(x_n)}{f'(x_n)}, \qquad n=0,1,\ldots,
\]
to approximate roots (or zeroes) of real-valued function $f(x)=0$ (cf. \cite{atkinson}).
Here the initial value $x_0$ is describe as the initial guess or iterate of the iteration process. 
Other applications in related sciences, such as in Theoretical Biology, also appear in countless literatures (see, e.g., \cite{cushing, din, elaydi} and other related literatures cited therein).
In a modeling setting, the two-dimensional \emph{competitive system} of nonlinear difference equations
	\[
		x_{n+1}=\frac{x_{n}}{a+y_{n}}, \quad y_{n+1}=\frac{y_{n}}{b+x_{n}}, \qquad n=0,1,\ldots,
	\]
represents the rule by which two discrete competitive populations reproduce from one generation to the next. 
In this context, the phase variables $x_n$ and $y_n$ denote population sizes during the $n$-th generation of a given population. 
The sequence or orbit $\{(x_n,y_n)\}_{n\in \mathbb{N}_0}$, where $\mathbb{N}_0:=\mathbb{N}\cup \{0\}$, describes how the populations evolve over time. 
The system is describe as a competition between the populations because of the fact that the transition function for each population is a decreasing function of the other population size. 
In \cite{hassel}, Hassell and Comins studied a discrete (difference) single age-class model for two-species competition and investigated the stability properties of its solutions.
Evidently, difference equations appear abundantly in nature and have many applications in biology, ecology and epidemiology \cite{kapur}, economy \cite{neusser}, physics (see, e.g., \cite{courant}) and so on.
Because of its importance in science, these equations will no doubt attract the attentions of many scientists and especially the interest of mathematicians.
Another reason which possibly creates greater interest in this topic is the challenge of understanding the behavior of solutions of these types of equations.
In fact, some difference equations appear simple in forms but the behavior of their solutions, as we have said, are difficult and sometimes too complex to completely undrestand.

In this study, we aim to describe the qualitative behavior of solutions of some systems of nonlinear difference equations.
This investigation is motivated by earlier results in this area, wherein different systems of difference equations have been solved in closed-forms.
We have seen that several earlier results concerning the form of solutions of these types of equations were established through the principle of induction,
and no further justifications on how these formulae were obtained  are stated.
Our specific purpose here is to fill these gaps and provide theoretical explanations on how these results can be established analytically.
We mention that we did a similar investigation in \cite{rabago} wherein we studied the system of nonlinear difference equations
\begin{equation}
\label{eqrabago}
x_n=\frac{a}{y_{n-p}},\qquad y_n=\frac{by_{n-p}}{x_{n-q}y_{n-q}},\qquad n=0,1,\ldots.
\end{equation}
where $q$ is a positive integer with $p < q$, $p \nmid q$, and $p$ is an even number, both $a$ and $b$ are nonzero real constants, and
the initial values $x_{-q+1},x_{-q+2},\ldots,x_0$, $y_{-q+1},y_{-q+2},\ldots,y_0$ are nonzero real numbers.
In our previous work \cite{rabago}, we are able to answer completely an open question raised by Yang et al. in \cite{ycs}.
That is, we are able to describe the behavior of solutions of the system \eqref{eqrabago} when $b=a$, $p < q$, $p \nmid q$ and $p$ is an even number.
We have particularly found that every solution of this system when $b=a$, with $p>0$, is unbounded whenever $q$ is odd.
However, a periodic solution of the given system occurs when $\gcd(p,q)>1$ and $p/\gcd(p,q)$ is odd. 
In this case, the period of the solution appears to be equal to the least common multiple of $p$ and $2q$.
On the other hand, a similar behavior as for the case when $q$ is odd was observed when $\gcd(p,q)>1$ and $p/\gcd(p,q)$ is even. 
It is worth mentioning that other related papers also appeared in literature before ours \cite{rabago}, 
see, for example, Cinar \cite{c1}, Cinar and Yal\c{c}inkaya \cite{cy1,cy2,cy3}, Cinar et al. \cite{cyr1}, Elsayed \cite{e1}, 
Iri\^{c}anin and Liu \cite{bratislav}, \"{O}zban \cite{oz1,oz2,oz3} and Yang et al. \cite{ylb1}.
With our set objective, we shall revisit two recent studies (see \cite{bayram} and \cite{elsayed}) and present more results regarding the topic discussed in these papers.
In the first part of our paper, we consider the system of difference equations
\begin{equation}
\label{system}
x_{n+1} =\frac{1}{y_{n-k}}, \quad y_{n+1}=\frac{x_{n-k}}{y_{n-k}}, \qquad n=0,1,\ldots\tag{Sys. 1}
\end{equation}
with positive initial values $\{x_i\}_{i=-k}^0$ and $\{y_i\}_{i=-k}^0$.
The above system was first studied in \cite{bayram} by Bayram and Da\d{s}, 
wherein they obtained a result regarding the behavior of solution of system \ref{system}:
\begin{theorem}[cf. \cite{bayram}]
\label{main}
Let $\{x_i\}_{i=-k}^0$ and $\{y_i\}_{i=-k}^0$ be positive real numbers for some integer $k\geq1$. 
Then, every solution of the system \ref{system} is periodic with period $3(k + 1)$.
\end{theorem}
The above result was established in \cite{bayram} in a quite laborious manner.
More specifically, the authors \cite{bayram} were able to show the validity of the above theorem by computing iteratively the form of solutions of the system \ref{system}.
Evidently, this method requires much time and effort before one can able to find exactly the periodicity of solution of the given system. 
One of our object here in this work is to provide an alternative approach in determining the periodicity of periodic solution of system \ref{system}.
This approach, as we will see later on, shall be more appreciated because of its simplicity and elegance.

Another recent work which shall be revisited here is due to Elsayed \cite{elsayed}.    
In \cite{elsayed}, Elsayed investigated the behavior of solutions of several instances of the following system of nonlinear difference equations
\[
x_{n+1} = \frac{1}{y_{n-k}}, \quad y_{n+1}=\frac{y_{n-k}}{x_{n-m}y_{n-m}}, \qquad n=0,1,\ldots.
\]
He particularly studied the solutions of the above system with (i) $(k,m)=(2r,0)$, (ii) $(k,m)=(2r-1,0)$ and (iii) $(k,m)=(2r-1,1)$ where $r$ is some positive integer. 
It was shown that the system, with the given cases, have periodic and unbounded solutions depending on some set requirements on the initial values.
Most, but not all, of these results was substantiated by Elsayed with the same approach used by Bayram and Da\d{s} \cite{bayram}.
However, we have found that there is an error in one of the result stated in \cite{elsayed}.
Perhaps because the result was not supported by any analytical justifications. 
Hence, in this work, we shall correct this mistake and offer some theoretical explanations on the results presented in \cite{elsayed}. 

The rest of the paper is organized as follows: in the next section (Section 2), we shall provide a theoretical explanation of Theorem \ref{main}, which can also be viewed as another way to prove the theorem.
In Section 3, we investigate the behavior of solutions of the system of difference equations
\[
x_{n+1} = \frac{1}{y_{n-k}}, \quad y_{n+1}=\frac{y_{n-k}}{x_{n-m}y_{n-m}}, \qquad n=0,1,\ldots,
\]
and determine the conditions for which the system will have a periodic or unbounded solutions. 
In Section 4, we establish, in a way different from Elsayed \cite{elsayed}, the form of solutions of the above system when $k$ is odd and $m=0,1$ (see systems \ref{s3} and \ref{s4}).
For completeness, we shall accompany our results with several illustrations.
Finally, in Section 5, we provide a summary of our present work. 

Before we further proceed in our work, we note that term sequences and difference equations will be interchangeably use throughout our discussion.
\section{Proof of Theorem \ref{main}}
Throughout the proof we denote the solution of system \ref{system} with $\{(x_n,y_n)\}:=\{(x_n,y_n)\}_{n=-k}^{\infty}$.
To begin with, we eliminate $x_{n+1}$ in system \ref{system} to obtain the one-dimensional difference equation 
\begin{equation*}
y_{n+1}=\frac{1}{y_{n-k}y_{n-2k-1}}\qquad\Longleftrightarrow\qquad y_{n+1}y_{n-k}y_{n-2k-1} =1, \qquad n=0,1,\ldots
\end{equation*}
Given that the initial values $\{y_i\}_{i=-k}^0$ are positive real numbers, we can then take the natural logarithm of both sides of the latter equation yielding the relation
\begin{equation}
\label{log}
 \ln y_{n+1} + \ln y_{n-k} + \ln y_{n-2k-1} =0, \qquad n=0,1,\ldots. 
\end{equation}
Defining $a_n:=\ln y_n$, we can equivalently write \eqref{log} as follows
\begin{equation}
\label{rec}
a_{n+1}+a_{n-k}+a_{n-2k-1}=0,  \qquad n=0,1,\ldots. 
\end{equation}
The above relation is obviously a linear recurrence equation of degree $2k+2$.
Using the ansatz $a_n=\lambda^n$ for some $\lambda \in \mathbb{C}$, equation \eqref{rec} therefore has the characteristic equation
\[
P(\lambda):=\lambda^{2k+2}+\lambda^{k+1}+1=0.
\]
Letting $\mu=\lambda^{k+1}$ in $P(\lambda)=0$, we obtain the transformed equation
\[
Q(\mu)=\mu^2+\mu+1.
\]
Now, it is evident that the complex conjugate roots of $Q(\mu)=0$ is of order $3$ (i.e., if $\mu^{\ast}$ satisfies $Q(\mu^{\ast})=0$, then $(\mu^{\ast})^3=1$), all of which are found on the unit disk $|z|\leq 1$.
This can be seen easily from the fact that the polynomial equation $(\mu-1)Q(\mu)=\mu^3-1=0$ has roots $\mu=\exp(2l\pi i/3), l=0,1,2$.
It follows that 
\[
\lambda_l=\left\{\exp\left(\frac{2l\pi i}{3}\right)\right\}^{1/(k+1)}, \qquad l=0,1,2, 
\]
are the roots of $P(\lambda)=0$. 
Clearly, since the roots of $Q(\mu)=0$ are distinct, then so are the roots of $P(\lambda)=0$, i.e., all roots of $P(\lambda)=0$ are simple.
Hence, the sequence $\{a_n\}$ has the explicit formula of the form
\begin{equation}
\label{aform}
a_n = c_1 \lambda_1^n + c_2 \lambda_2^n + \ldots +\ldots +c_{k+1} \lambda_{k+1}^n+\ldots + c_d \lambda_d^n, \qquad n=0,1,\ldots,
\end{equation}
for some real numbers $\{c_j\}_{j=1}^d$, where $\{\lambda_j\}_{j=1}^d$, $d:=3(k+1)$, are roots of the polynomial equation $P(\lambda)=0$.
With this notation, we see that $\lambda_j^d=1$, for all $j=1,2,\ldots,d$. 
In fact, we have $\lambda_j^{dn+t}=\lambda_j^t$, for all $n=0,1,\ldots$ and $t\in\{0,1,\ldots,d-1\}$.
So, from formula \eqref{aform}, we have
\[
a_{dn+t} = \sum_{j=1}^d c_j \lambda_j^{dn+t}=\sum_{j=1}^d c_j \lambda_j^t = a_t, \qquad n=0,1,\ldots,
\]
for all $t \in\{0,1,\ldots,d-1\}$.
Thus, the sequence $\{a_n\}$ is periodic of period $d$.
Going back to the relation $a_n=\ln y_n$, we see that the sequence $\{y_n\}$ is periodic with period $3(k+1)$.
Similarly, $\{x_n\}$ is periodic with the same period.
This proves Theorem \ref{main}. 

\begin{remark}
We remark that the proof of Theorem \ref{main} also provides a method for determining the existence of periodic solution of similar system of difference equations reducible to linear recurrence equations.
Once that the corresponding characteristic equation of the transformed equation contains a repeated root, then we can argue that every solution of the given system is unbounded.
This idea shall be later on use to deal with the behavior of solutions of another system of nonlinear difference equations.   
\end{remark}


\section{More on the system $x_{n+1} =1/y_{n-k}, \ y_{n+1}=y_{n-k}/(x_{n-m}y_{n-m})$}

Let $k$ and $m$ be fixed positive integers. In this section we consider the system of difference equations 
\begin{equation}
\label{system2}
x_{n+1} = \frac{1}{y_{n-k}}, \quad y_{n+1}=\frac{y_{n-k}}{x_{n-m}y_{n-m}}, \qquad n=0,1,\ldots\tag{Sys. 2}
\end{equation}
with nonzero real initial conditions $\{x_n\}_{n=-\nu}^0$ and $\{y_n\}_{n=-\nu}^0$, where $\nu=\max\{k,m\}$.
Hereafter, we denote the solution of the above system by $\{(x_n,y_n)\}:=\{(x_n,y_n)\}_{m=-\nu}^{\infty}$.
Some particular cases of \ref{system2} were first studied in \cite{elsayed}:
\begin{align}
x_{n+1} &=\frac{1}{y_{n-2r}}, \quad y_{n+1}=\frac{y_{n-2r}}{x_{n}y_{n}}, \qquad n=0,1,\ldots,\label{s1}\tag{Sys. 2.1}\\
x_{n+1} &=\frac{1}{y_{n-2r}}, \quad y_{n+1}=\frac{y_{n-2r}}{x_{n-m}y_{n-m}}, \qquad n=0,1,\ldots,\label{s2}\tag{Sys. 2.2}\\
x_{n+1} &=\frac{1}{y_{n-2r-1}}, \quad y_{n+1}=\frac{y_{n-2r-1}}{x_{n}y_{n}}, \qquad n=0,1,\ldots,\label{s3}\tag{Sys. 2.3}\\
x_{n+1} &=\frac{1}{y_{n-2r-1}}, \quad y_{n+1}=\frac{y_{n-2r-1}}{x_{n-1}y_{n-1}}, \qquad n=0,1,\ldots,\label{s4}\tag{Sys. 2.4}
\end{align}
where $r$ are some positive integer.
It was shown by Elsayed in \cite{elsayed} that all solutions to system \ref{s1} are periodic with period $4r+2$, 
while he stated without proof that all solutions to system \ref{s2} are periodic with period $(2p+2)(2r+1)$ when $2p\neq r$ and periodic with period $(4r+2)$ when $2p=r$. 
Elsayed also showed in \cite{elsayed} that when $x_0y_0=1$, then all solutions of system \ref{s3} are periodic with period $(2r + 2)$.
Meanwhile, if $x_0y_0\neq 1$ in system \ref{s3}, then the solutions are unbounded.
This latter result is seen to be a consequence of the form of solutions of system \ref{s3} established in the paper.
Finally, he had shown in \cite{elsayed} that every solution to system \ref{s4} are also periodic of period $2r+2$ when $x_0y_0=x_{-1}y_{-1}=1$ and are unbounded when either $x_0y_0\neq 1$ or $x_{-1}y_{-1}\neq1$. 
The latter statement of this last-mentioned result, however, is erroneous, but can be fixed with additional conditions on $k$, as we shall see later on in our forthcoming discussion.

We emphasize that the results regarding the periodicity of periodic solutions of systems \ref{s1}--\ref{s4} were all established in \cite{elsayed} through the same approach used by Bayram and Da\d{s} \cite{bayram}. 
That is, by iteratively computing for the value of each $x_n$ (resp. $y_n$) until they arrive at a relation showing that $x_n=x_{n+p}$ (resp. $y_n=y_{n+p}$) for some integer $p>0$. 

In the sequel we shall determine the periodicity of periodic solutions of systems \ref{s1}--\ref{s4} through the same method we used in proving Theorem \ref{main}.
We accompany each of our results with illustrative examples where the initial values are randomly taken from the unit interval $(0,1)$.

Now, since systems \ref{s1}--\ref{s4} are particular cases of \ref{system2}, we go directly on the problem of determining the periodicity of periodic solution of system \ref{system2}, 
and then branch out from this system by treating individually all of its possible cases which will naturally arise from the original system.

\subsection{Case $m=k$} 
\begin{figure}
   \centering
    \scalebox{.35}{\includegraphics{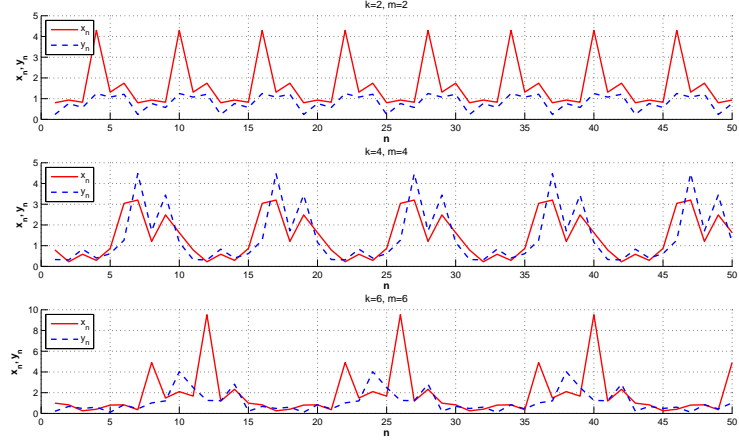}}
    \caption{The above figures illustrate the behavior of solutions of system \ref{system2} when $m=k$. 
    In these examples, we have considered the values $m=k\in\{2,4,6\}$.
    Notice that the illustrated behavior of solutions are periodic with period $6, 10$ and $14$, respectively for each cases $m=k=2,3$ and $4$.}\label{fig1}
\end{figure}

First, we suppose that $m=k$. 
Hence, we have the system
\[
x_{n+1} =\frac{1}{y_{n-k}}, \quad y_{n+1}=\frac{1}{x_{n-k}}, \qquad n=0,1,\ldots
\]
from which follows that
\[
y_{n+1}=y_{n-2k-1},  \qquad n=0,1,\ldots.
\]
This simply shows that if $m=k$, then every solution to \ref{system2} is periodic with period $2k+2$.
Figure \ref{fig1} illustrates the behavior of solutions of \ref{system2} when $m=k$.
\subsection{Case $m\neq k$} 
Suppose now that $m\neq k$ and $x_0y_0=1$. 
Then, multiplying $x_{n+1}$ and $y_{n+1}$, we get
\[
x_{n+1}y_{n+1} = \left( \frac{1}{y_{n-k}} \right)\left( \frac{y_{n-k}}{x_{n-m}y_{n-m}}\right) = \frac{1}{x_{n-m}y_{n-m}}, \qquad n=0,1,\ldots.
\]
If, in particular, $m=0$, then the above relation reduces to 
\[
x_{n+1}y_{n+1} = \frac{1}{x_ny_n}, \qquad n=0,1,\ldots.
\]
Given that $x_0y_0=1$, then it follows immediately that $x_ny_n=1$ for all $n\in \mathbb{N}_0$.
Going back to system \ref{system2}, we see that when $m=0$, then 
\[
y_{n+1}=y_{n-k}, \qquad n=0,1,\ldots,
 \] 
 for any integer $k>0$. 
This relation implies that the sequence $\{y_n\}$ is periodic with period $k+1$, which also means that $\{x_n\}$ is periodic with the same period.
Thus, every solution to system \ref{system2} when $m=0$ and $x_0y_0 =1$ is periodic with period $k+1$.

\begin{remark}
From our previous discussion, we can recover Elsayed's result for the behavior of solution of system \ref{system2} when $m=0$ and $k=2r+1$ for some integer $r \geq 0$.
That is, we see that all solutions of system \ref{s3} when $x_0y_0=1$ are indeed periodic with period $2r+2$ (cf. \cite[Theorem 2]{elsayed}).
\end{remark}

\begin{remark}
In general, if $m>0$ and $x_iy_i=1$ for all $i=-m,-m+1,\ldots,-1,0$, then system \ref{system2} is periodic with period $k+1$ (cf. \cite[Theorem 2, Proposition 2--(i)]{elsayed}).
\end{remark}

Now consider the case when $x_iy_i\neq 1$ for at least one $i\in\{-m,-m+1,\ldots,-1,0\}$, then in reference to system \ref{system2}, we obtain the relation
\[
y_{n+1}y_{n-m}=y_{n-k}y_{n-m-k-1}, \qquad n=0,1,\ldots
\]
which can be transformed into 
\begin{equation}
\label{an}
a_{n+1}+a_{n-m}-a_{n-k}-a_{n-m-k-1}=0, \qquad n=0,1,\ldots
\end{equation}
using the equation $a_n:=\ln y_n$.
The above recurrence equation has the characteristic equation 
\[
f(\lambda):
=(\lambda^{k+1}-1)(\lambda^{m+1}+1)=0,
\]
whose roots are given by 
\[
\lambda_j:=\exp\left\{\frac{2j\pi i}{k+1}\right\}, \ j=0,1,\ldots,k,
\quad\text{and}\quad
\hat{\lambda}_l:=\exp\left\{\frac{(2l+1)\pi i}{m+1}\right\}, \ l=0,1,\ldots,m.
\]
Evidently, all roots of $f(\lambda)=0$ will be simple if the inequality
\begin{equation}
\label{inequality}
2j(m+1)\neq(2l+1)(k+1)
 \end{equation}
holds for each $j,l\in \mathbb{N}_0$.
Henceforth, we consider several cases.\\

\begin{figure}[h!]
   \centering
    \scalebox{.35}{\includegraphics{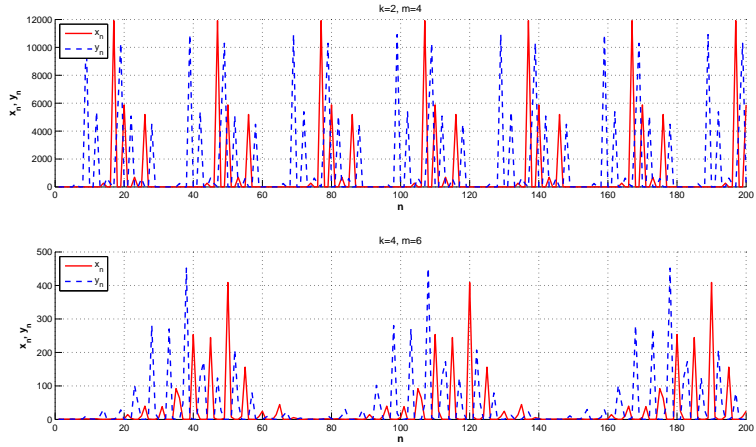}}
    \caption{The above figures illustrate the behavior of solutions of system \ref{system2} when $(m,k)=(2,4)$ (upper plot) and $(m,k)=(4,6)$ (lower plot), respectively.
    Observe in these two plots that the illustrated behavior of solutions of these examples are periodic with period $30$ and $70$, respectively.}\label{fig2}
\end{figure}

\begin{figure}[h!]
   \centering
    \scalebox{.35}{\includegraphics{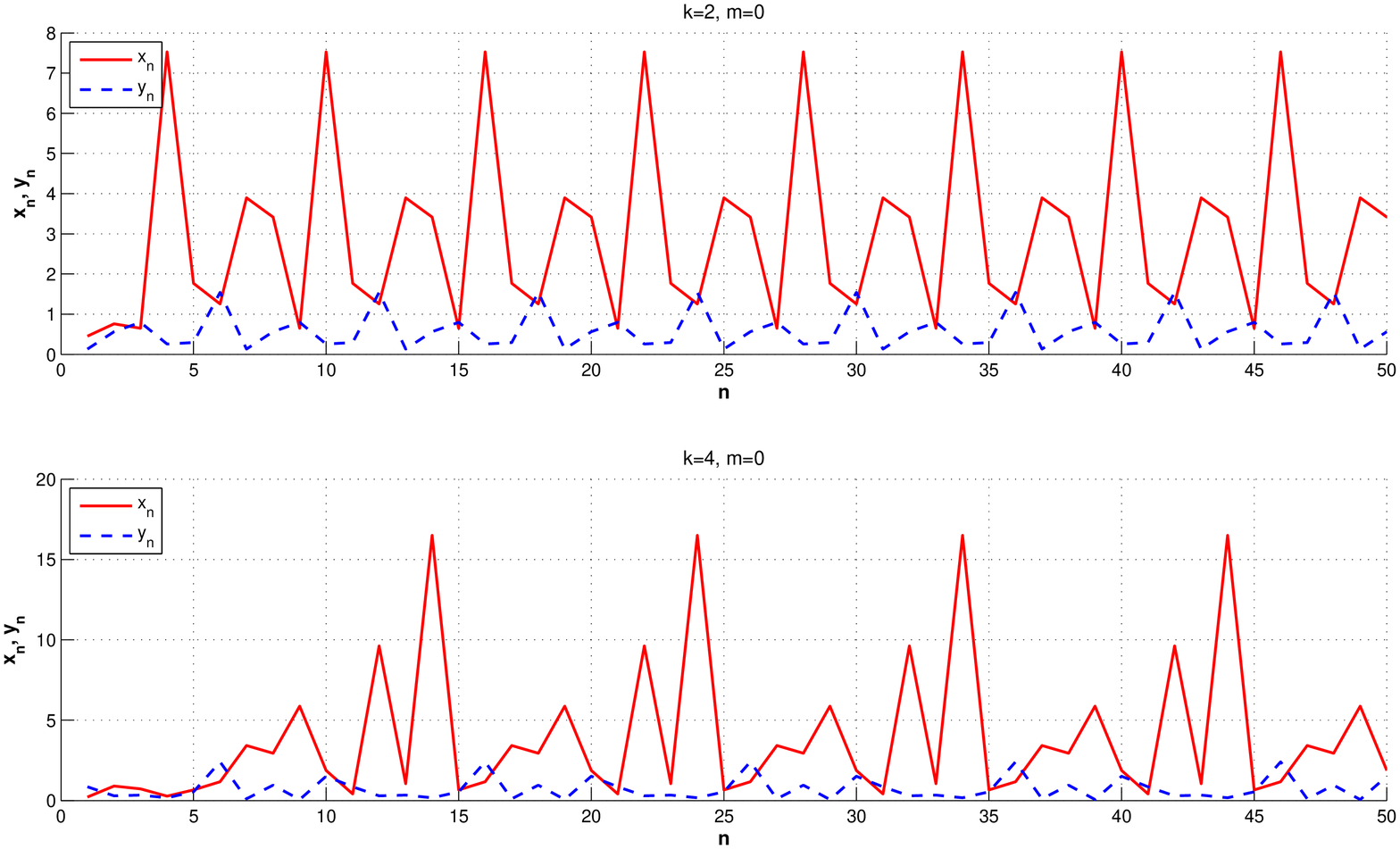}}
    \caption{The above figures illustrate the behavior of solutions of system \ref{system2} when $m=0$ for $k=2,4$, (respectively, upper and lower plot).
    In these cases, the solutions are periodic with period $6$ and $10$, respectively.}\label{fig3}
\end{figure}

\begin{figure}[h!]
   \centering
    \scalebox{.35}{\includegraphics{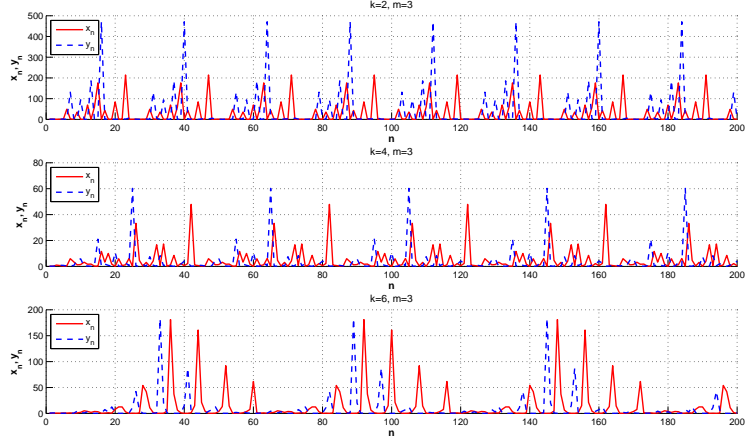}}
    \caption{The above figures illustrate the behavior of solutions of system \ref{system2} when $m$ is odd (in particular, $m=3$) for $k=2,4,6$, (respectively, upper, middle and lower plot).
    Evidently, the solutions are periodic with period $2(k+1)(m+1)=6(k+1)=24,30,56$, respectively for $k=2,4,6$.}\label{fig4}
\end{figure}

\underline{CASE 1}: Suppose $k$ is even. 
Then, obviously, inequality \eqref{inequality} holds true for each $j,l \in \mathbb{N}_0$. 
Hence, all roots of \ref{system} are simple, and the $n$-th term $a_n$ of the sequence satisfying the recurrence relation \ref{an} has the explicit formula of the form
\[
a_n=\sum_{j=1}^{k+1} c_j \lambda_j^n + \sum_{l=1}^{m+1}\hat{c}_l \hat{\lambda}_l^n, \qquad n=0,1,\ldots,
 \]   
 for some real numbers $c_1,c_2,\ldots,c_{k+1},\hat{c}_1,\hat{c}_2,\ldots$ and $\hat{c}_{m+1}$.
 Since all roots of $f(\lambda)=0$ lie on the unit disk $|z|\leq 1$, and $\lambda_j^{k+1}=1$ for each $j=0,1,\ldots, k$ and $\hat{\lambda}_l^{m+1}=-1$ for each $l=0,1,\ldots, m$,
 then arguing as in the proof of Theorem \ref{main}, we get
\[
a_{2(k+1)(m+1)n+t}=a_t,  \qquad n=0,1,\ldots,
\]
for all $t\in\{0,1,\ldots, 2(k+1)(m+1)-1\}$.
The above equation clearly indicates that when $k$ is even, then every solutions of the sequence $\{a_n\}$ is periodic with period $2(k+1)(m+1)$.
Thus, every solution $\{(x_n,y_n)\}$ of system \ref{system2} when $k$ is even is periodic with period $2(k+1)(m+1)$.
Some illustrations for this case are shown in Figures \ref{fig2}, \ref{fig3} and \ref{fig4}.

\begin{remark}
We emphasize that the result we obtain here in this case regarding the periodicity of periodic solution of system \ref{system2} agrees with \cite{elsayed}.
More precisely, we have validated in our previous discussion Elsayed's claim in \cite[Proposition]{elsayed}.
We also remark that, as an immediate consequence, every positive solution of system \ref{s1} is periodic with period $4r+2$ (cf. \cite[Theorem 1]{elsayed}).
\end{remark}

\underline{CASE 2}: Suppose $k$ is odd and $m$ is even, i.e., $k=2r+1$ and $m=2s$ for some integers $r,s>0$. 
Then, from inequality \eqref{inequality}, we must have
\[
2j(2s+1) \neq (2l+1)(2r+2)\quad\Longleftrightarrow\quad j(2s+1)\neq(2l+1)(r+1) 
\]
for each $j,l\in \mathbb{N}_0$, where $r$ and $s$ are fixed positive integers, for $f(\lambda)=0$ to have simple roots.
However, the equation $j(2s+1)=(2l+1)(r+1)$ always holds by choosing $j=r+1$ and $l=s$.
This implies that $f(\lambda)=0$ always has a root of order two. 
Without loss of generality (WLOG), suppose $\lambda_1$ is a root of order two and let $\hat{\lambda}_1=\lambda_1$.
Then, the explicit formula for $a_n$ takes the form
\[
a_n= c_1 \lambda_1^n + \hat{c}_1 n \hat{\lambda}_1^n + \sum_{j=2}^{k+1} c_j \lambda_j^n + \sum_{l=2}^{m+1}\hat{c}_l \hat{\lambda}_l^n, \qquad n=0,1,\ldots.
\]
In order to show that every solution of the system is non-periodic (i.e., has an unbounded solution), it suffices to prove that it has a subsequence that tends to  infinity (or perhaps, converges to zero).
Since $k$ and $m$ are of different parity, then $\gcd(k,m)=1$ which in turn implies that
\[
\lambda_j^{(k+1)(m+1)}=\hat{\lambda}_l^{(k+1)(m+1)}=1,
\]
 for all $j=1,2,\ldots,k+1$ and $l=1,2,\ldots,m+1$, respectively. 
 Moreover, we have
\[
a_{(k+1)(m+1)n}= c_1 + \hat{c}_1 (k+1)(m+1)n + \sum_{j=2}^{k+1} c_j + \sum_{l=2}^{m+1}\hat{c}_l, \qquad n=0,1,\ldots.
\]
Suppose (WLOG) that $\hat{c}>0$, then it follows that $a_{2(k+1)(m+1)n}$ is increasing for all $n \geq N$ for some sufficiently large $N\in \mathbb{N}_0$.
In fact, for all $n \geq N$, we'll have
\[
a_{(k+1)(m+1)n} \longrightarrow \infty \qquad\text{as}\qquad n \longrightarrow \infty.
\]
Thus, the subsequence $\{a_{(k+1)(m+1)n}\}$ tends to infinity as $n$ increases.
Therefore, the sequence $\{a_n\}$ is unbounded.
Going back to the relation $a_n=\ln y_n$, we conclude that there is a subsequence of $\{y_n\}$ (resp. $\{x_n\}$) which tends to infinity (resp. converges to zero) exponentially.
This implies that every solution $\{(x_n,y_n)\}$ of system \ref{system2} is unbounded, and therefore nonperiodic when $k$ is odd and $m$ is even. 
\begin{figure}[h!]
   \centering
    \scalebox{.35}{\includegraphics{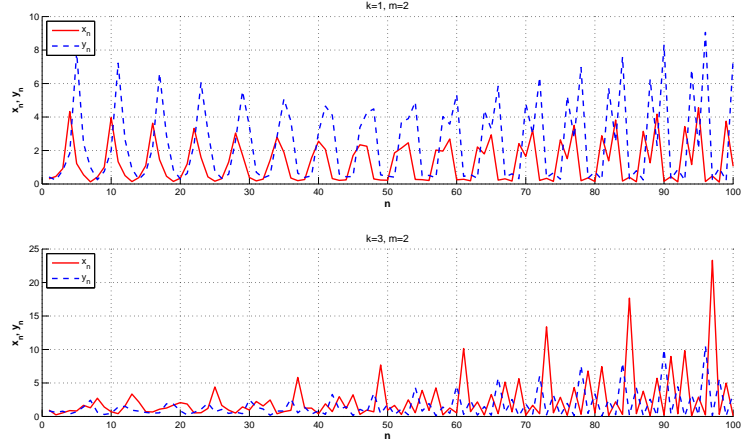}}\label{fig5}
    \caption{The above figures illustrate the behavior of solutions of system \ref{system2} when $k$ is odd and $m$ is even.
    Observe from the upper plot that the illustrated solution of the system starts to increase after $N=50$, while it is already evident in the lower plot that the particular solution of the system is exponentially increasing.}
\end{figure}


\begin{remark}
We mention that Elsayed's result in \cite[Theorem 2--(ii)]{elsayed} is a consequence of our previous result.
In particular, when $x_0y_0\neq 1$, and $k=2r+1$ and $m=0$ (which is exactly system \ref{s3}), then all solutions of system \ref{system2} are unbounded. 
\end{remark}

\underline{CASE 3}: Suppose $k$ and $m$ are both odd, i.e., $k=2r+1$ and $m=2s+1$ for some integers $r,s>0$. 
Hence, we require
\[
2j(2s+2) \neq (2l+1)(2r+2)\quad\Longleftrightarrow\quad 2j(s+1)\neq(2l+1)(r+1) 
\]
to hold for each $j,l\in \mathbb{N}_0$, where $r$ and $s$ are fixed positive integers, so that $f(\lambda)=0$ will have simple roots.
However, the equation $2j(s+1)=(2l+1)(r+1)$ may hold true by choosing $j=(r+1)/2$ and $l=s/2$, or perhaps, when $2(s+1)=r+1$.
The former statement is only possible when $r$ is odd and $s$ is even.
So, in contrary, if $r$ is even or $s$ is odd, and $2(s+1)\neq r+1$, then all roots of the equation $f(\lambda)=0$ are simple.
Arguing as in the previous case, we conclude that if $\{(x_n,y_n)\}$ is a solution to system \ref{system2} with $k$ and $m$ both odd, then one the following statements are true:
\begin{enumerate}
\item[(i)] if $k \equiv -1\ ({\rm mod}\ 4)$ and $m \equiv 1\ ({\rm mod}\ 4)$, or $2(m+1)=k+1$, then $\{(x_n,y_n)\}$ is unbounded, i.e. system \ref{system2} will never have a periodic solution.
\item[(ii)] if $k \equiv 1\ ({\rm mod}\ 4)$ or $m \equiv -1\ ({\rm mod}\ 4)$, and $2(m+1)\neq k+1$, then $\{(x_n,y_n)\}$ is periodic with period equal to the least common multiple of $k+1$ and $2(m+1)$.
\end{enumerate} 
The periodicity we mention in statement (ii) follows from the fact that $\hat{\lambda}_l^{m+1}=-1$ for all $l=0,1,\ldots,m$ and hence, $\hat{\lambda}_l^{2(m+1)}=1$.
However, since $k+1$ is even, then it suffices to take the least common multiple of $k+1$ and $2(m+1)$, say $\rho$, to have $\hat{\lambda}_l^{\rho}=1$ for each $l=0,1,\ldots,m$.
\begin{figure}[h!]
   \centering
    \scalebox{.35}{\includegraphics{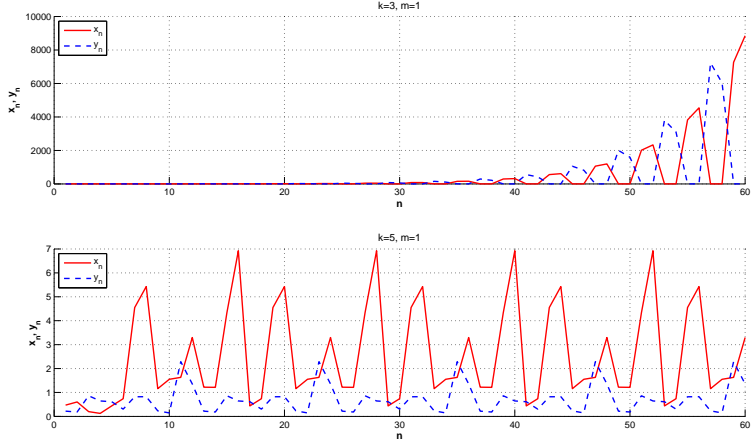}}\label{fig7}
    \caption{The above figures illustrate the behavior of solutions of system \ref{system2} when $k$ and $m$ are both odd. 
    Referring to the upper plot, we see that when $k=3 \equiv -1 \ ({\rm mod}\ 4)$ and $m \equiv 1\ ({\rm mod}\ 4)$ we'll have an unbounded solution.
    This agrees with the conclusion in Case 3--(i).
    On the other hand, it is evident that $k=5$ and $m=1$ satisfies the conditions in Case 3--(ii), so we see that the solution illustrated on the lower plot is periodic with period $12$ which exactly equates to ${\rm lcm(6,2)}$, 
    where ${\rm lcm}(\cdot,\cdot)$ denotes the least common multiple of two numbers.}
\end{figure}

\begin{remark}
Elsayed stated in \cite[Proposition 2--(ii)]{elsayed} that every solution of system \ref{system2} when $k$ is odd and $m=1$ (cf. system \ref{s4}) is unbounded.
However, a counterexample (refer to {\rm Figure  (6)}, lower plot) would show that this statement is not true in general. 
\end{remark}

We now summarize our previous discussion in the following theorem.
\begin{theorem}
\label{summary}
Let $m$ and $k$ be any positive integers and consider the system of nonlinear difference equations \ref{system}, with positive initial conditions $\{x_n\}_{n=-\nu}^0$ and $\{y_n\}_{n=-\nu}^0$, where $\nu:=\max\{k,m\}$.
Then, a solution $\{(x_n,y_n)\}_{n=-\nu}^{\infty}$ to \ref{system2} behaves accordingly as follows:

\begin{enumerate}
\item[{\rm (i)}] if $m=k$, then $\{(x_n,y_n)\}_{n=-\nu}^{\infty}$ is periodic with period $2k+2$.

\item[{\rm(ii)}] If $m\neq k$ and $x_iy_i =1$ for all $i=-m,-m+1,\ldots, -1,0$, then $\{(x_n,y_n)\}_{n=-\nu}^{\infty}$ is periodic with period $k+1$.

\item[{\rm(iii)}] If $m\neq k, x_iy_i\neq 1$ for at least one $i\in\{-m,-m+1,\ldots,-1,0\}$, and 
   \begin{enumerate}
	 \item[{\rm(a)}] $k$ is even, then $\{(x_n,y_n)\}_{n=-\nu}^{\infty}$ is periodic with period $2(k+1)(m+1)$.
	 \item[{\rm(b)}] $k$ is odd while $m$ is even, then $\{(x_n,y_n)\}_{n=-\nu}^{\infty}$ is unbounded.
	 \item[{\rm(c)}] $k \equiv -1\ ({\rm mod}\ 4)$ and $m \equiv 1\ ({\rm mod}\ 4)$, or $2(m+1)=k+1$, then $\{(x_n,y_n)\}_{n=-\nu}^{\infty}$ is unbounded.
	 \item[{\rm(d)}] $k \equiv 1\ ({\rm mod}\ 4)$ or $m \equiv -1\ ({\rm mod}\ 4)$, and $2(m+1)\neq k+1$, 
	then $\{(x_n,y_n)\}_{n=-\nu}^{\infty}$ is periodic with period equal to the least common multiple of $k+1$ and $2(m+1)$.
  \end{enumerate}
\end{enumerate}
\end{theorem}

\section{Forms of Solutions of Systems \ref{s3} and \ref{s4}}
Here we find a closed-form solution of systems \ref{s3} and \ref{s4}.
We established these expressions through an analytical approach and not with the usual induction method, thus providing a theoretical explanation on Elsayed's result presented in \cite[Theorem 2--(ii) and Proposition 2--(ii)]{elsayed}.
We only consider the case when all initial values are positive real number. 
The method, however, can be followed inductively to extend the results to the more general case.

To have an idea about the form of solution of systems \ref{s3} and \ref{s4}, we state in advance the following results: 

\begin{theorem}[\cite{elsayed}]
\label{t2}
Let $x_0, y_{-(2r+1)},\ldots,y_{-1}$, and $y_0$ be nonzero real numbers such that $A:=x_0y_0 \neq 1$. 
Then, every solution of system \ref{s3} takes the form
\[
x_{(2r+2)n+s} = 
\begin{cases}
\dfrac{A^n}{y_{-(2r+2)+s}}, & \text{when} \ s\ \text{is odd},\\[1em]
\dfrac{1}{A^ny_{-(2r+2)+s}} ,& \text{when} \ s\ \text{is even},
\end{cases}
\quad s=1,2,\ldots,2r+2,
\]
and
\[
y_{(2r+2)n+s} = 
\begin{cases}
\dfrac{y_s}{A^n}, & \text{when} \ s\ \text{is odd},\\[1em]
A^n y_s ,& \text{when} \ s\ \text{is even},
\end{cases}
\quad s=-2r-1,-2r,\ldots,0,
\]
for all $n=0,1,\ldots$.
\end{theorem}


\begin{theorem}[\cite{elsayed}]
\label{t3}
Let $x_{-1}, x_0, y_{-(2r+1)},\ldots,y_{-1}$, and $y_0$ be nonzero real numbers such that $A:=x_0y_0 \neq 1$ or $B:=x_{-1}y_{-1}\neq1$. 
Suppose that the conditions in (iii)-(d) are satisfied.
Then, the solution of system \ref{s4} takes the form
\[
x_{(2r+2)n+s} = 
\begin{cases}
\dfrac{B^{n\sin\{s\pi/2\}}}{y_{-(2r+2)+s}}, & \text{when} \ s\ \text{is odd},\\[1em]
\dfrac{A^{n\cos\{(s+2)\pi/2\}}}{y_{-(2r+2)+s}} ,& \text{when} \ s\ \text{is even},
\end{cases}
\quad s=1,2,\ldots,2r+2,
\]
and
\[
y_{(2r+2)n+s} = 
\begin{cases}
\dfrac{y_s}{B^{n\sin\{s\pi/2\}}}, & \text{when} \ s\ \text{is odd},\\[1em]
\dfrac{y_s}{A^{n\cos\{(s+2)\pi/2\}}} ,& \text{when} \ s\ \text{is even},
\end{cases}
\quad s=-2r-1,-2r,\ldots,0,
\]
for all $n=0,1,\ldots$.
\end{theorem}

\begin{remark}
As we have mentioned before, Theorem \ref{t2} has already been proven in \cite{elsayed} through induction method.
On the other hand, Theorem \ref{t3} can be compared with \cite[Proposition 2]{elsayed} 
but we emphasize that the conclusion that every solution when $k$ is odd (with no further conditions on $k$) is unbounded is not correct.
However, this was easily fixed in Theorem \ref{summary} by imposing additional conditions on $k$ and $m$.
The solution, as we saw from (iii)-(d) in Theorem \ref{summary}, will only be periodic if $k\neq m$, $k \equiv 1\ ({\rm mod}\ 4)$ or $m \equiv -1\ ({\rm mod}\ 4)$, and $2(m+1)\neq k+1$.
\end{remark}

Consider the product 
\begin{equation}
\label{xy}
x_{n+1}y_{n+1} = \frac{1}{x_{n-m}y_{n-m}}, \qquad n=0,1,\ldots.
\end{equation}
Observe that if $x_jy_j$ is positive for each $j=-m,-m+1,\ldots,-1,0$, then so is $x_ny_n$ for all $n\in\mathbb{N}$.
Hence, we may take the logarithm of both sides of equation \eqref{xy}. 
That is, upon letting $a_n:=\ln x_ny_n$, we obtain the one-dimensional linear difference equation
\[
a_{n+1} + a_{n-m} = 0, \qquad n=0,1,\ldots.
\]
Clearly, using the ansatz $a_n=\lambda^n$, the above recurrence relation has the characteristic equation 
\[
g(\lambda):=\lambda^{m+1}+1=0.
\]

Now in the succeeding discussion we offer a theoretical approach to Theorems \ref{t2} and \ref{t3}.\\


\begin{proof}[Proof to Theorem \ref{t2}] 
If $m=0$, then we obtain the polynomial equation $g(\lambda)=\lambda+1=0$.
Hence, $\lambda=-1$. 
So the sequence $\{a_n\}$ has the explicit formula given by
\[
a_n= {\rm const.} \times (-1)^n, \qquad n=0,1,\ldots,
\]
where the constant must be equal to $a_0$.
Since $\ln x_ny_n=a_n$, then we have
\[
x_ny_n = \left\{ x_0 y_0\right\}^{ (-1)^n}, \qquad n=0,1,\ldots.
\]
But, $x_n=1/y_{n-k-1}$, so 
\[
y_n=\frac{\left\{ x_0 y_0\right\}^{ (-1)^n}}{x_n} = \left\{ x_0 y_0\right\}^{ (-1)^n} y_{n-k-1}, \qquad n=0,1,\ldots.
\]
Replacing $n$ by $(k+1)n+s$ in the above equation, we get
\[
y_{(k+1)n+s}= \left\{ x_0 y_0\right\}^{ (-1)^{(k+1)n+s}} y_{(k+1)(n-1)+s}, \qquad n,s=0,1,\ldots.
\]
Iterating the left hand side of the above equation and then substituting $k$ by $2r+1$, we obtain
\[
y_{(2r+2)n+s}= \left\{ x_0 y_0\right\}^{ (-1)^sn} y_s, \qquad n,s=0,1,\ldots.
\]
Evidently, the above form of solution for $y_{(2r+2)n+s}$ can be equivalently express as
\[
y_{(2r+2)n+s} = 
\begin{cases}
\dfrac{y_s}{A^n}, & \text{when} \ s\ \text{is odd},\\[1em]
A^n y_s ,& \text{when} \ s\ \text{is even},
\end{cases}
\quad s=-2r-1,-2r,\ldots,-1,0.
\]
The desired form, on the other hand, for $x_{(2r+2)n+s}$ follows from the relation $x_n=\{x_0y_0\}^{(-1)^sn}/y_n$.\\
This proves Theorem \ref{t2}.  
\end{proof}

Now we turn on the proof of Theorem \ref{t3}.
In this part, we take into account our results in Theorem \ref{summary}--(iii). 
Particularly, the result we stated in case (c) and (d).
In the proof, we assume that $k\equiv -1\ ({\rm mod}\ 4)$ so that the first condition in (iii)-(c) of Theorem \ref{summary} is satisfied.\\

\begin{proof}[Proof to Theorem \ref{t3}]
If $m=1$, then $g(\lambda)=\lambda^2+1=0$ which has two complex conjugate roots $\lambda_{\pm}=\pm i$.
Since the two roots are evidently distinct, then $a_n$ can be express as
\[
a_n = c_1 i^n + c_2 (-i)^n, \qquad n=0,1,\ldots, 
 \] 
 where $(c_1,c_2)$ is the solution of the system $c_1+c_2=a_0$ and $c_1i-c_2i=a_1$, 
 i.e., 
 \[
 c_1=\frac{1}{2}(a_0- i a_1)
 \qquad \text{and}\qquad
  c_2=\frac{1}{2}(a_0+i a_1).
 \]
 Hence,
 \begin{align*}
 a_n &= \frac{1}{2}\left\{ (a_0- i a_1) i^n + (a_0+i a_1) (-i)^n\right\}\\
  &= a_0\frac{i^n+(-i)^n}{2} + a_1\frac{i^{n-1} + (-i)^{n-1}}{2}, \qquad n=0,1,\ldots.
  \end{align*}
Notice that 
\[
\xi(n):=\frac{i^n+(-i)^n}{2} = 
\begin{cases}
1, & \text{if} \ n\equiv0\ (\text{mod}\ 4),\\[0.5em]
0, & \text{if} \ n\equiv1\ (\text{mod}\ 4),\\[0.5em]
-1, & \text{if} \ n\equiv2\ (\text{mod}\ 4),\\[0.5em]
0, & \text{if} \ n\equiv3\ (\text{mod}\ 4),
\end{cases}
\quad n=0,1,\ldots.
\]
Hence, 
\[
a_n = a_0 \xi(n) + a_1\xi(n-1), \qquad n=0,1,\ldots.
\]
Now, since $a_0=\ln x_0 y_0=:\ln A$ and $a_1=\ln x_1y_1=-\ln x_{-1}y_{-1}=:-\ln B$, then we have
\[
x_ny_n =A ^{\xi(n)} B^{-\xi(n-1)}, \qquad n=0,1,\ldots. 
\]
Because $x_n=1/y_{n-k-1}$, then the above relation is equivalent to
\begin{equation}
\label{yform}
y_n =A^{\xi(n)} B^{-\xi(n-1)}y_{n-k-1}, \qquad n=0,1,\ldots. 
\end{equation}
Replacing $n$ by $(2r+2)n+s$, we get
\begin{equation}
\label{yns}
y_{(2r+2)n+s} = A^{\xi((2r+2)n+s)} B^{-\xi((2r+2)n+s-1)} y_{(2r+2)(n-1)+s}, \qquad n=0,1,\ldots.
\end{equation}

We consider two possibilities for $n$ as follows:\\

\underline{POSSIBILITY 1}: If $n$ is even, say $n=2N$ for some $N \in\mathbb{N}_0$, then $\xi(4(r+1)N+s) =\xi(s)$. 
So, iterating the left hand side (LHS) of \eqref{yns} with $n=2N$, we get
\[
y_{(2r+2)n+s} =  A^{n\xi(s)} B^{-n\xi(s-1)} y_s.
\]
Hence, if $n=2N$, we have the formula
\[
y_{(2r+2)n+s} = 
\begin{cases}
y_sB^{-n\xi(s-1)}, & \text{when} \ s\ \text{is odd},\\[1em]
A^{n\xi(s)} y_s ,& \text{when} \ s\ \text{is even},
\end{cases}
\quad s=-2r-1,-2r,\ldots,-1,0.
\]
Here follows that
\[
x_{(2r+2)n+s} = 
\begin{cases}
\dfrac{B^{n\xi(s-1)}}{y_{-(2r+2)+s}}, & \text{when} \ s\ \text{is odd},\\[1em]
\dfrac{A^{-n\xi(s)}}{y_{-(2r+2)+s}} ,& \text{when} \ s\ \text{is even},
\end{cases}
\quad s=1,2,\ldots, 2r+2.
\]
It can be verified easily by induction that the following identities hold:
\[
\xi(s) = -\cos\left\{ \frac{(s+2)\pi}{2}\right\}
\quad\text{and}\quad 
\xi(s-1) = \sin\left\{ \frac{s\pi}{2}\right\}
\]
for all $s=-2r-1,-2r,\ldots,-1,0$.
Using these identities, we get the desired form for $y_{(2r+2)n+s} $ and $x_{(2r+2)n+s}$ as in Theorem \ref{t3} with $n$ even.\\

\underline{POSSIBILITY 2}: Now suppose that $n$ is odd, i.e., $n=2N+1$, then 
\[
\xi((2r+2)(2N+1)+s) =\xi(2r+2+s).
\]
Iterating the LHS of equation \eqref{yns} as we did in the previous case, we obtain
\begin{equation}
\label{ynrs}
y_{(2r+2)n+s} =  A^{n\xi(2r+2+s)} B^{-n\xi(2r+1+s)} y_s.
\end{equation}
Replacing $s$ by $-(2r+2)+s$, we get
\begin{equation}
\label{yform1}
y_{(2r+2)(n-1)+s} =  A^{n\xi(s)} B^{-n\xi(s-1)} y_{-(2r+2)+s}.
\end{equation}
However, from equation \eqref{ynrs}, letting $n=1$ and replacing $s$ by $-(2r+2)+s$, we'll have the relation
\[
y_s=A^{\xi(s)} B^{-\xi(s-1)} y_{-(2r+2)+s}
\quad \Longleftrightarrow\quad
y_{-(2r+2)+s}=A^{-\xi(s)} B^{\xi(s-1)}y_s. 
\]
Substituting this latter expression for $y_{-(2r+2)+s}$ in equation \eqref{yform1}, we obtain
\begin{align*}
y_{(2r+2)(n-1)+s} 
&=  A^{n\xi(s)} B^{-n\xi(s-1)} A^{\xi(s)} B^{-\xi(s-1)}y_s\\
&=A^{(n-1)\xi(s)} B^{-(n-1)\xi(s-1)} y_s, \qquad n=1,2,\ldots
\end{align*}
Observe that $n-1=(2N+1)-1=2N$, i.e., $n-1$ is even.
Thus, going back to the previous possibility considered, we get the desired result.
This proves the second case, completing the proof of Theorem \ref{t3}. 
\end{proof}

In what follows, we established the form of solutions of system \ref{s4} wherein the conditions in (iii)-(d) of Theorem \ref{summary} are satisfied,
i.e., we find the closed-form solution of the system
\begin{equation}
x_{n+1} =\frac{1}{y_{n-k}}, \quad y_{n+1}=\frac{y_{n-k}}{x_{n-1}y_{n-1}}, \qquad n=0,1,\ldots,\label{s5}\tag{Sys. 2.5},
\end{equation}
with $k \equiv 1\ ({\rm mod}\ 4)$ (with $k\neq m=1$). 

So, suppose $k \equiv 1\ ({\rm mod}\ 4)$ (with $k\neq m=1$). 
Then, obviously $2(m+1) \neq k+1$ (since $m=1$ and $k+1 =4k'+2$ for some $k' \in \mathbb{N}$).
Note that the conditions we imposed on $k$ agree with that of (iii)-(d) of Theorem \ref{summary}.
Hence, every solution $\{(x_n,y_n)\}_{n=\nu}^{\infty}$ of system \ref{s5} 
with $m=1$ and $k \equiv 1\ ({\rm mod}\ 4)$ is periodic with period equal to ${\rm lcm}(k+1,2(m+1))=(4k'+2,4)=4(2k'+1)$.
Having this idea in mind, we can explicitly obtain the solutions $\{(x_n,y_n)\}$ in closed-form.
To do this, we use the fact that the solution $\{(x_n,y_n)\}$ are periodic with period $4(2k'+1)$ and that $\max\{k,m\}=4k'+1$
so as to obtain the relations:
\[
y_{n+(4k'+2)}=y_{n-(4k'+2)}\quad \text{and}\quad x_{n+(4k'+2)}=x_{n-(4k'+2)}, \qquad \text{for all} \ n=0,1,2,\ldots.
\]
However, with reference to system \ref{s5}, we have 
\[
x_n=\frac{1}{y_{n-(4k'+2)}}=y_{n+(4k'+2)}\quad\text{and}\quad y_n=\frac{1}{x_{n+(4k'+2)}}=\frac{1}{x_{n-(4k'+2)}}, \qquad n=0,1,\ldots.
\]
Once again, from equation \ref{s5}, we have
  \[
  y_{n+(4k'+2)}=\frac{y_n}{x_{n+(4k'+2)-2}y_{n+(4k'+2)-2}}, \qquad n=0,1,2,\ldots.
  \]
But, in view of the equation $x_{n+1}y_{n+1} = 1/(x_ny_n)$, we have the relation
\[
x_ny_n=x_{n-2}y_{n-2}, \qquad n=0,1,2,\ldots,
\]
which means that the product sequence $\{x_ny_n\}_{n=\nu}^{\infty}$ is periodic of period $2$.
Thus, the closed-form solutions of system \ref{s5} with $k \equiv 1 \ ({\rm mod}\ 4)$ can now be formally stated as follows.
\begin{theorem}
\label{t5}
Let $x_{-1}, x_0, y_{-(2r+1)},\ldots,y_{-1}$, and $y_0$ be nonzero real numbers such that $A:=x_0y_0 \neq 1$ or $B:=x_{-1}y_{-1}\neq1$.
Furthermore, let $k=4k'+1$ for some $k'\in\mathbb{N}$ and suppose that the conditions in (iii)-(d) are satisfied. 
Also, defined $\Omega_1:=\{1,2,\ldots,4k'+2\}$ and $\Omega_2:=\{4k'+3,4k'+4,\ldots,\rho:=4(2k'+1)\}$. 
Then, every solution of system \ref{s5} takes the form
\[
x_n =
\begin{cases}
\dfrac{1}{y_{-(4k'+2)+s}},				&\text{if} \ n \ ({\rm mod}\ \rho) \in \Omega_1,\\[1em] 
\dfrac{A^{\xi(s)}B^{-\xi(s-1)}}{y_{-(4k'+2)+s}},	&\text{if} \ n \ ({\rm mod}\ \rho) \in \Omega_2, 
\end{cases}
\]
and
\[
y_n =
\begin{cases}
\dfrac{y_{-(4k'+2)+s}}{A^{-\xi(s)}B^{\xi(s-1)}},	&\text{if} \ n \ ({\rm mod}\ \rho) \in \Omega_1,\\[1em] 
y_{-(4k'+2)+s},						&\text{if} \ n \ ({\rm mod}\ \rho) \in \Omega_2, 
\end{cases}
\]
where $s\equiv n \ ({\rm mod}\ 4(2k'+1))$, for all $n=1,2,\ldots$.
\end{theorem}

To illustrate our previous result, we consider the system 
\[
x_{n+1} =\frac{1}{y_{n-5}}, \quad y_{n+1}=\frac{y_{n-5}}{x_{n-1}y_{n-1}}, \qquad n=0,1,\ldots,\label{s6}\tag{Sys. 2.6},
\]
with real initial values $x_{-5}, x_{-4}, x_{-3}, x_{-2}, x_{-1}, x_0, y_{-5}, y_{-4}, y_{-3}, y_{-2}, y_{-1}$ and $y_0$.
So, in view of Theorem \ref{t5}, the solutions of system \ref{s5} are given as follows:
\[
x_n =
\begin{cases}
\dfrac{1}{y_{-6+s}},								&\text{if} \ n \ ({\rm mod}\ 12) \in \{1,2,3,4,5,6\},\\[1em] 
\dfrac{\{x_0y_0\}^{\xi(s)}\{x_{-1}y_{-1}\}^{-\xi(s-1)}}{y_{-6+s}},	&\text{if} \ n \ ({\rm mod}\ 12) \in \{7,8,9,10,11,0\}, 
\end{cases}
\]
and 
\[
y_n =
\begin{cases}
\dfrac{y_{-6+s}}{\{x_0y_0\}^{-\xi(s)}\{x_{-1}y_{-1}\}^{\xi(s-1)}},	&\text{if} \ n \ ({\rm mod}\ 12) \in \{1,2,3,4,5,6\},\\[1em] 
y_{-6+s},										&\text{if} \ n \ ({\rm mod}\ 12) \in \{7,8,9,10,11,0\}, 
\end{cases}
\]
where $s\equiv n \ ({\rm mod}\ 12)$, for all $n=1,2,\ldots$. 
That is, every solution of system \ref{s6} takes the following form:
\[
\{x_n\}_{n=1}^{\infty}=\left\{\frac{1}{y_{-5}}, \frac{1}{y_{-4}}, \frac{1}{y_{-3}}, \frac{1}{y_{-2}}, \frac{1}{y_{-1}}, \frac{1}{y_0},
\frac{x_{-1}y_{-1}}{y_{-5}}, \frac{x_0y_0}{y_{-4}}, \frac{1}{x_{-1}y_{-1} y_{-3}}, \frac{1}{x_0y_0 y_{-2}}, x_{-1}, x_0, \ldots\right\}
\]
and
\[
\{y_n\}_{n=1}^{\infty}=\left\{\frac{y_{-5}}{x_{-1}y_{-1}}, \frac{y_{-4}}{x_0y_0}, x_{-1}y_{-1} y_{-3}, x_0y_0 y_{-2}, \frac{1}{x_{-1}}, \frac{1}{x_0},
y_{-5}, y_{-4}, y_{-3}, y_{-2}, y_{-1}, y_0, \ldots\right\}.
\]

The next two results are immediate consequences of Theorems \ref{t2} and \ref{t3}, respectively.

\begin{corollary}
Let $x_0, y_{-(2r+1)},\ldots,y_{-1}$, and $y_0$ be nonzero real numbers such that $A:=x_0y_0 \neq 1$. 
Then, every solution of system \ref{s3} are unbounded.
\end{corollary}

{\bf Proof}
Suppose (WLOG) that $A>1$. 
Then, from Theorem \ref{t2}, we have
\[
x_{(2r+2)n+s} \longrightarrow 
\begin{cases}
\infty,& \text{when} \ s\ \text{is odd},\\[0.75em]
0,& \text{when} \ s\ \text{is even},
\end{cases}
\quad s=1,2,\ldots,2r+2,
\]
and on the other hand, we have
\[
y_{(2r+2)n+s} \longrightarrow  
\begin{cases}
0, & \text{when} \ s\ \text{is odd},\\[0.75em]
\infty,& \text{when} \ s\ \text{is even},
\end{cases}
\quad s=-2r-1,-2r,\ldots,0,
\]
from which it is evident that the solution $\{(x_n,y_n)\}$ of system \ref{s3} is unbounded.

\begin{corollary}
Let $x_{-1}, x_0, y_{-(2r+1)},\ldots,y_{-1}$ and $y_0$ be nonzero real numbers such that $A:=x_0y_0 \neq 1$ or $B:=x_{-1}y_{-1}\neq1$.
Suppose that the conditions in (iii)-(c) are satisfied.
Then, every solution of system \ref{s4} is unbounded.
\end{corollary}

We end our paper with a summary of most important results found here in this investigation.
\section{Summary}

We have demonstrated analytically that every positive solution of the system 
\[
x_{n+1} =\frac{1}{y_{n-k}}, \quad y_{n+1}=\frac{x_{n-k}}{y_{n-k}}, \qquad n=0,1,\ldots
\]
is periodic with period $3(k+1)$. 
Our result validated, in an alternative approach, the outcome issued by Bayram and Da\d{s}' (2010) in \cite{bayram}.
We have also established, in an elegant fashion, the behavior of solutions of the system of difference equations
\[
x_{n+1} =\frac{1}{y_{n-k}}, \quad y_{n+1}=\frac{y_{n-k}}{x_{n-m}y_{n-m}}, \qquad n=0,1,\ldots
\]
with nonzero real initial conditions $\{x_n\}_{n=-\nu}^0$ and $\{y_n\}_{n=-\nu}^0$, where $\nu=\max\{k,m\}$.
As special cases of the latter system, we have confirmed, all except for one result, those that were presented in \cite{elsayed}.
Particularly, in contrast to Elsayed's claim in \cite[Proposition 2-(ii)]{elsayed}, we have substantiated in this work that the given system will have a periodic solution 
when $k$ is odd and $m=1$ and if the conditions $k \equiv 1\ ({\rm mod}\ 4)$ or $m \equiv -1\ ({\rm mod}\ 4)$, and $2(m+1)\neq k+1$ are satisfied.
In this instance, we have shown that the periodic solution has periodicity equal to the least common multiple of $k+1$ and $m+1$.
The system, as we have justified and theoretically explained, may also have periodic solutions in the following cases: 
\begin{enumerate}
\item[(i)] when $m=k$; 
\item[(ii)] when $m\neq k$ and $x_iy_i =1$ for all $i=-m,-m+1,\ldots, -1,0$;
\item[(iii)] when $m\neq k$ and $k$ is even with $x_iy_i\neq 1$ for at least one $i\in\{-m,-m+1,\ldots,-1,0\}$.
\end{enumerate}
In these situations, the solutions are shown to be periodic with periodicity taking value $2k+2$, $k+1$ and $2(k+1)(m+1)$, respectively. 
On the other hand, the system exhibits unbounded solutions (iv) when $k$ is odd while $m$ is even and $x_iy_i\neq 1$, for at least one $i\in\{-m,-m+1,\ldots,-1,0\}$,
and (v) when $m\neq k$ with $k \equiv -1\ ({\rm mod}\ 4)$ and $m \equiv 1\ ({\rm mod}\ 4)$, or $2(m+1)=k+1$.
Consequently, more interesting results can be established for other related systems following the idea presented in this work.
Thus, in our next study, we continue investigating the behavior of solutions of other related systems of nonlinear difference equations and contribute more in this developing topic.



\end{document}